# An Algorithm for Computing Ideals and Conjugacy Classes of Subalgebras of Borel Subalgebras

Nimra Sher Asghar[a*], Hassan Azad[b]

[a]*School of Natural Sciences, National University of Sciences and Technology, H-12, 44000, Islamabad, Pakistan*
[b]*Abdus Salam School of Mathematical Sciences, GCU, 54600, Lahore, Pakistan*

**Abstract**

In this article, we present a constructive procedure for determining all ideals of the Borel subalgebra of a complex semisimple Lie algebra from its root system or, equivalently, its Dynkin diagram. The proposed algorithmic approach has been implemented in Maple. The Maple code has been tested on Borel subalgebras associated with both classical and exceptional Lie algebras. An interactive procedure is also given to determine conjugacy classes of all subalgebras in the derived algebra of the Borel subalgebra.

## 1. Introduction

The main aim of this paper is to develop a constructive procedure for identifying all ideals of the Borel subalgebra of a given complex semisimple Lie algebra from its root system and implement this algorithmic approach in Maple. This is used to determine the conjugacy classes of subalgebras of dimensions at most 3 in derived algebras of Borel subalgebras in rank-2 semisimple algebras. The restriction to rank-2 algebras is made because of its applications to the classifications of all subalgebras of vector fields in three variables that have a proper Levi decomposition.

Although we can compute all ideals in classical algebras of a given rank and in exceptional Lie algebras, the resulting lists are too long to be presented in a paper. However, these can be obtained by substituting the type of the root system in the given code. As an illustration of this, we have listed all the abelian ideals of type $F_4$ algebras.

*Corresponding author
 *Email addresses:* `nimra.phdmaths22sns@student.nust.edu.pk` (Nimra Sher Asghar[a]),
`hassan.azad@sms.edu.pk` (Hassan Azad[b])

Earlier work on this topic is due to several authors. Notably, Peterson [1], Cellini and Papi [2], Rö̈hrle [3, 5], Suter [4], and Panyushev [5–8]. Additionally, a different approach to finding conjugacy classes is given by Amata and Oliveri [9]. In view of the importance of this subject, it is desirable to give alternative approaches.

This paper is organized as follows: In Section 2, we present a constructive procedure for determining all ideals of the Borel subalgebra of a complex semisimple Lie algebra. Section 3 details the implementation of this procedure in Maple. Section 4 demonstrates the execution of the developed code using rank-2 semisimple Lie algebras as examples. Additionally, we utilize the computed ideals of the derived algebra of Borel subalgebras to determine the conjugacy classes in these algebras. Lastly, in Section 5, we list all the abelian ideals of type F4 algebras. We refer the reader to [10] for all results used from the theory of root systems.

## 2. Algorithmic Approach to Ideals of Borel Subalgebra

This section presents an algorithmic method to systematically determine all ideals of a Borel subalgebra within a complex simple Lie algebra, using Proposition 2.1 of [11].

Let R be the root system of a simple Lie algebra, $R^+$ be a system of positive roots, and S the corresponding system of simple roots. The main result of this section is the following:

**Proposition 2.1.**

1. *Let*

$$\mathfrak{g} = T \oplus \sum_{\alpha \in R^+} \mathbb{C} \cdot X_\alpha$$

    *be the Borel subalgebra of a complex semisimple Lie algebra. Then, every 1-dimensional ideal of $\mathfrak{g}$ is $\langle X_h \rangle$, where $[X_a, X_h] = 0$ for all simple roots a. In particular, if R is irreducible, then $\mathfrak{g}$ has a unique one-dimensional ideal, given by $\mathbb{C} \cdot X_h$, where h is the highest root.*

2. *Let $J \subsetneq \mathfrak{g}'$ be a d-dimensional ideal of $\mathfrak{g}$ properly contained in $\mathfrak{g}'$. Let $R_J$ be the roots of $T$ in $J$. Then, a 1-dimensional ideal of $\mathfrak{g}/J$ in $\mathfrak{g}'/J$ is given by $\langle \overline{X_\alpha} \rangle$, with $\alpha \in R^+/R_J$ with $[X_a, X_\alpha] = 0$ or $a + \alpha \in R_J$ for all simple roots a.*

*Proof.*

1. From the decomposition

$$\mathfrak{g} = T \oplus \sum_{\alpha \in R^+} \mathbb{C} \cdot X_\alpha,$$

we see that $\mathfrak{r}$ is self centralizing and therefore $Z(\mathfrak{g}) = 0$. Thus a 1-dimensional ideal of $\mathfrak{g}$ is contained in

$$\mathfrak{g}' = \sum_{\alpha \in R^+} \mathbb{C} \cdot X_\alpha.$$

As every ideal of $\mathfrak{g}$ is $\mathfrak{r}$-invariant, we see that a 1-dimensional ideal of $\mathfrak{g}$ is of the form $\mathbb{C} \cdot X_\alpha$ with

$$[X_r, X_\alpha] = 0,$$

for all $r \in R^+$. This is so if and only if $[X_a, X_\alpha] = 0$ for all simple roots $a$.

Therefore, in case R is irreducible,

$$\langle X_\alpha \rangle = \langle X_h \rangle,$$

where, $h$ is the highest root.

2. Let $\mathbf{J} \subset \mathfrak{g}'$ be a $d$-dimensional ideal of $\mathfrak{g}$, assume that $\mathbf{J} /= \mathfrak{g}'$. We want to determine all 1-dimensional ideals of $\mathfrak{g}/\mathbf{J}$ contained in $\mathfrak{g}'/\mathbf{J}$.

As J is $\mathfrak{r}$-invariant, it is a sum of root spaces, say

$$\mathbf{J} = \sum_{\alpha \in R_J} \mathfrak{g}_\alpha,$$

and therefore

$$\mathfrak{g}'/\mathbf{J} = \sum_{\alpha \in R^+ \setminus R_J} \bar{\mathfrak{g}}_\alpha \quad \text{and} \quad \mathfrak{g}/\mathbf{J} = \bar{\mathfrak{r}} + \sum_{\alpha \in R^+ \setminus R_J} \bar{\mathfrak{g}}_\alpha.$$

Suppose $\bar{t}_0 + \sum_{\alpha \in R^+ \setminus R_J} c_\alpha \bar{X}_\alpha$ is in $Z(\mathfrak{g}/\mathbf{J})$. Thus for all $t \in \mathfrak{r}$,

$$\sum_{\alpha \in R^+ \setminus R_J} c_\alpha \cdot \alpha(t) \cdot \bar{X}_\alpha = 0.$$

But then $c_\alpha = 0$ for all $\alpha \in R^+ \setminus R_J$. This means that

$$Z(\mathfrak{g}/\mathbf{J}) \subset \bar{\mathfrak{r}}.$$

Therefore, $Z(\mathfrak{g}/\mathbf{J})^\mathsf{T} \mathfrak{g}'/\mathbf{J} = 0$. Thus a 1-dimensional ideal of $\mathfrak{g}/\mathbf{J}$ in $\mathfrak{g}'/\mathbf{J}$ is given by an element of $X$ in

$$Z(\mathfrak{g}'/\mathbf{J}) = Z((\mathfrak{g}/\mathbf{J})'),$$

which is a common eigenvector of $\bar{\mathfrak{g}}/\langle \bar{\mathfrak{g}}', Z(\bar{\mathfrak{g}}) \rangle$. Moreover,

$$Z(\mathbf{g}/\mathbf{J}) = \{\bar{t} \in \bar{\tau} \mid \alpha(t) = 0 \,\forall\, \alpha \in R^+ \setminus R_J \}.$$

Now $Z(\mathbf{g}')$ is invariant under $\bar{\tau}$ and thus it is spanned by $\bar{X}_\alpha$ with $\alpha \in R^+ \setminus R_J$ and $[X_\alpha, X_a] = 0$ or $\alpha + a \in R_J$ for all simple roots $a$.

Thus in case R is irreducible, starting with the ideal $\langle X_h \rangle$, where $h$ is the highest root, we get all ideals of $\mathbf{g}$ in $\mathbf{g}'$.

$\square$

**Corollary 2.2.** *All ideals J on $\mathbf{g}$ are of the form*

$$\mathbf{J} = S + \mathbf{J} \cap \mathbf{g}',$$

*where S is a subalgebra of*

$$\bigcap_{\beta \in R^+ \setminus R_J} kernel(\beta).$$

*with $R_J$ being the roots of $\tau$ in J.*

For constructing examples, we recall for the convenience of the reader a construction procedure to obtain all positive roots starting from the simple roots

$$S = \{\alpha_1, \ldots, \alpha_l\}.$$

Compute for each pair of simple roots $(\alpha, \beta)$, with $\langle \alpha, \beta \rangle < 0$, the roots $\sigma_\alpha(\beta)$, where $\sigma_\alpha$ is the reflection given by $\alpha$. Remove repetitions and call this set $S_1$. For each pair $(r, \alpha)$, $r \in S_1, \alpha \in S$ with $\langle r, \alpha \rangle < 0$, compute

$$\sigma_\alpha(r) = r - \langle r, \alpha \rangle \alpha.$$

Remove repetitions and continue till a set $S_k$ is obtained for which

$$\langle r, \alpha \rangle \geq 0,$$

for all $r \in S_k$ and $\alpha \in S$. At this stage, all the roots have been obtained.

## 3. Computational Aspects of Ideals in Borel Subalgebra

This section presents a Maple implementation based on Sec. 2, providing an algorithmic method to systematically determine all ideals of a given Borel subalgebra in a complex semisimple Lie algebra. The Maple packages `DifferentialGeometry`, `LieAlgebras`, `Library`, and `plots` are used, and they are loaded using the `with` command. All the functions used in the code are defined in the appendix. The developed code is as follows:

```
Type_of_Root := " INPUT ":   # Input a string, a root type "A", "B", "C", "D", "E", "F", "G"
n := " INPUT ": #positive integer n without string
Dynkin_Diagram := DynkinDiagram(Type_of_Root, n):
display(Dynkin_Diagram);
Root_Type := parse(Type_of_Root)[n];
Cartan_Matrix := CartanMatrix(Type_of_Root, n):
Positive_Roots := PositiveRoots(Cartan_Matrix):
Simple_Roots := SimpleRoots(Positive_Roots):
R_oot_s := AlphaForge(Positive_Roots, Simple_Roots):
Positive_Roots := R_oot_s[1];
Simple_Roots := R_oot_s[2];
Cartan_Subalgebra := ConstructCartanSubalgebra(Positive_Roots, Simple_Roots):
Derived_Algebra := ConstructDerivedAlgebra(Positive_Roots, Simple_Roots);
Borel_Subalgebra := ConstructBorelSubalgebra(Cartan_Subalgebra, Derived_Algebra);
AllIdeals := []:
One_Dimensional_Ideals := DetermineOneDimensionalIdeal(Positive_Roots, Simple_Roots):
Current_Dimensional_Ideals := One_Dimensional_Ideals:
AllIdeals := {op(AllIdeals), op(One_Dimensional_Ideals)}:
    while Current_Dimensional_Ideals <> {} do
        NextIdeal := []:
        for element in Current_Dimensional_Ideals do
            Positive_roots := `minus`({op(Positive_Roots)}, {op(RootsExtractor(element))});
            Positive_roots := [op(Positive_roots)];
            OtherIdeal:= DetermineOtherIdeals(element, Positive_roots, Simple_Roots);
            NextIdeal := {op(NextIdeal), op(OtherIdeal)};
            NextIdeal := map(sortList, NextIdeal):
        end do:
        Current_Dimensional_Ideals := convert(NextIdeal, set);
        Current_Dimensional_Ideals := {op(Current_Dimensional_Ideals)};
        AllIdeals := {op(AllIdeals), op(Current_Dimensional_Ideals)}
    end do:
```

```
32  Ideals_of_DA_of_Borel_Subalgebra := AllIdeals;
33  basislistCartan := BasisListCartan(Cartan_Subalgebra):
34  cefflist := CoefficientListCartan([seq(a[i], i = 1 .. nops(Cartan_Subalgebra))]):
35  Ideals_B_DA := CalculateDotProducts(cefflist, basislistCartan):
36  Ideals_containing_DA := []:
37      for elements in Ideals_B_DA do
38          result := [op(elements), op(Derived_Algebra)];
39          Ideals_containing_DA := {result, op(Ideals_containing_DA)}
40      end do:
41  Ideals_containing_DA:
42  All_Ideals_of_Borel_Subalgebra := {Borel_Subalgebra, op(Ideals_containing_DA), op(
        Ideals_of_DA_of_Borel_Subalgebra)};
```

Code 1: Ideals Determination of Borel Subalgebras

The provided code begins by taking the root type and rank ($n$) of the Lie algebra as inputs, followed by the construction of the Dynkin diagram, Cartan matrix, positive roots, and simple roots using appropriate Maple commands. To refine these sets of roots, the function `AlphaForge` is employed. Based on these structures, the Cartan subalgebra, the derived algebra $\mathfrak{g}'$, and the Borel subalgebra are constructed systematically.

To identify the one-dimensional ideals of the Borel subalgebra, the `DetermineOneDimensionalIdeal` function is constructed using the Proposition 2.1(1). This function takes two inputs: *proot*, representing the set of positive roots, and *sroot*, representing the set of simple roots. It systematically checks each root vector $X[i]$ in *proot* to verify whether it commutes with all simple root vectors $X[j] \in$ *sroot*, ensuring that it forms a one-dimensional ideal. The verification of the commutator condition is performed using the `lieB` function, which computes the Lie bracket by checking the sum of the roots associated with $X[i]$ and $X[j]$. If the sum is a valid root in *proot*, the bracket returns the corresponding root vector; otherwise, it returns zero. For each $X[i]$, if the condition holds for all $X[j]$, $X[i]$ is added to the list `OneDimensionalIdeals`. Finally, the function `DetermineOneDimensionalIdeal` returns the complete set of one-dimensional ideals.

The function `DetermineOtherIdeals` is designed to determine the higher-dimensional ideals of the Borel subalgebra, based on the concept of one-dimensional ideals. It takes three inputs: $I\_deal$, the current $d$-dimensional ideal; *proot*, the set of positive roots not yet included in the current ideal; and *sroot*, the set of simple roots. The function systematically checks which root vectors $X[k]$ from *proot* can be appended to $L\,deal$ to form the next ideal. It iterates over each root $X[k]$ in *proot* and checks the Lie bracket $[X[l], X[k]]$ for all $X[l] \in$ *sroot* using the func-tion `lieB`. If the bracket is zero for all $X[l]$, it implies that $X[k]$ commutes with the simple root vectors or satisfies the condition $a + \alpha \in \mathbf{R_J}$ from the Proposition 2.1(2). In this case, $X[k]$ is appended to $I\_deal$ to form a higher-dimensional ideal, which is then stored in `nextIdeal`. This process is integrated into an algorithm that iteratively constructs all possible ideals of the Borel subalgebra. The loop begins with the one-dimensional ideals, computed earlier using the `DetermineOneDimensionalIdeal` function, as `Current_Dimensional_Ideals`. At each iteration, `Current_Dimensional_Ideals` is updated by removing root vectors already included in the previously identified ideals, ensuring that each step progressively builds upon the prior structure. To construct higher-dimensional ideals, `DetermineOtherIdeals` function is used. The process continues to reach the derived algebra $\mathbf{g}'$.

Once the ideals of $\mathbf{g}'$ are obtained, the ideals of $\mathbf{g}$ containing $\mathbf{g}'$ are constructed. These ideals correspond to the subspaces of the abelian quotient $\mathbf{g}/\mathbf{g}'$. If $dim(\mathbf{g}/\mathbf{g}') > 1$, these are infinitely many ideals.

The iterative nature of the algorithm guarantees a comprehensive classification of all ideals in the Borel subalgebra, including those that contain the derived algebra. This systematic approach aligns with the theoretical framework, ensuring both completeness and accuracy in the identification of ideals.

## 4. Ideals and Conjugacy Classes in rank-2 Semisimple Lie Algebras

In this section, we present some relevant examples of complex semisimple Lie algebras and demonstrate how the Maple-based code allows us to determine all ideals of Borel subalgebras. We concentrate only on rank-2 as the results of higher ranks are too long to be given in a paper. Results for higher rank algebras can be obtained by substituting the type of the root system in the code.

### 4.1. *Borel subalgebras of type $A_2$*

We consider the $A_2$ algebra, which corresponds to $sl(3, \mathbb{C})$. The Maple code is executed with the root type `"A"` and $n = 2$ as inputs. The output, shown in Figure 1, provides the complete list of all ideals, illustrating the results obtained from the code.

○ ─────── ○
α₁         α₂

$$Root\_Type := A_2$$
$$Positive\_Roots := [\alpha 1, \alpha 2, \alpha 1 + \alpha 2]$$
$$Simple\_Roots := [\alpha 1, \alpha 2]$$
$$Derived\_Algebra := [X_{\alpha 1}, X_{\alpha 2}, X_{\alpha 1 + \alpha 2}]$$
$$Borel\_Subalgebra := [H_{\alpha 1}, H_{\alpha 2}, X_{\alpha 1}, X_{\alpha 2}, X_{\alpha 1 + \alpha 2}]$$
$$Ideals\_of\_DA\_of\_Borel\_Subalgebra := \{[X_{\alpha 1 + \alpha 2}], [X_{\alpha 1}, X_{\alpha 1 + \alpha 2}], [X_{\alpha 2}, X_{\alpha 1 + \alpha 2}], [X_{\alpha 1}, X_{\alpha 2}, X_{\alpha 1 + \alpha 2}]\}$$
$$All\_Ideals\_of\_Borel\_Subalgebra := \{[X_{\alpha 1 + \alpha 2}], [X_{\alpha 1}, X_{\alpha 1 + \alpha 2}], [X_{\alpha 2}, X_{\alpha 1 + \alpha 2}], [X_{\alpha 1}, X_{\alpha 2}, X_{\alpha 1 + \alpha 2}],$$
$$[H_{\alpha 2}, X_{\alpha 1}, X_{\alpha 2}, X_{\alpha 1 + \alpha 2}], [H_{\alpha 2} a_2 + H_{\alpha 1}, X_{\alpha 1}, X_{\alpha 2}, X_{\alpha 1 + \alpha 2}], [H_{\alpha 1}, H_{\alpha 2}, X_{\alpha 1}, X_{\alpha 2}, X_{\alpha 1 + \alpha 2}]\}$$

Figure 1: Ideals of a Borel subalgebras of type $A_2$

From figure 1, we observe the following steps:

- The Dynkin diagram for $A_2$ is plotted which consists of two nodes connected by a single edge, indicating a rank-2 algebra with two simple roots.

- Positive roots are calculated as $R^+ = \{\alpha 1, \alpha 2, \alpha 1 + \alpha 2\}$.

- Simple roots are identified as $S = \{\alpha 1, \alpha 2\}$.

- The Borel subalgebra is constructed as $\mathfrak{g} = \{H_{\alpha 1}, H_{\alpha 2}, X_{\alpha 1}, X_{\alpha 2}, X_{\alpha 1+\alpha 2}\}$.

- The derived algebra is determined as $\mathfrak{g}' = \{X_{\alpha 1}, X_{\alpha 2}, X_{\alpha 1+\alpha 2}\}$.

The code iteratively constructs and identifies all dimensional ideals, providing a complete classification of the ideals for $A_2$. Since $A_2$ is irreducible, its Borel subalgebra contains a unique one-dimensional ideal i.e.

$$\langle X_{\alpha 1+\alpha 2}\rangle.$$

The two-dimensional ideals of $\mathfrak{g}$ in $\mathfrak{g}'$ are

$$\langle X_{\alpha 1}, X_{\alpha 1+\alpha 2}\rangle \text{ and } \langle X_{\alpha 2}, X_{\alpha 1+\alpha 2}\rangle.$$

The three-dimensional ideal of $\mathfrak{g}$ in $\mathfrak{g}'$ is

$$\langle X_{\alpha_1}, X_{\alpha_2}, X_{\alpha_1+\alpha_2}\rangle = \mathfrak{g}'.$$

The code then identifies the ideals containing $\mathfrak{g}'$. To find these, we use the Cartan algebra $\mathfrak{g}/\mathfrak{g}'$, which is abelian with a basis $\langle H_{\alpha_1}, H_{\alpha_2}\rangle$. Thus, the four dimensional ideals containing $\mathfrak{g}'$ are of the form

$$\langle \lambda H_{\alpha_1} + \mu H_{\alpha_2}, X_{\alpha_1}, X_{\alpha_2}, X_{\alpha_1+\alpha_2}\rangle.$$

Specifically, these include $\langle H_{\alpha_1} + a_2 H_{\alpha_2}, X_{\alpha_1}, X_{\alpha_2}, X_{\alpha_1+\alpha_2}\rangle$ and $\langle H_{\alpha_2}, X_{\alpha_1}, X_{\alpha_2}, X_{\alpha_1+\alpha_2}\rangle$. This implies that there are infinitely many ideals in $\mathfrak{g}$. Additionally, the Borel subalgebra itself is also an ideal. By ideal we always mean a nontrivial ideal.

Extending the analysis of resultant ideals, we now aim to determine the conjugacy classes of subalgebras of a Borel subalgebra, focusing on its derived subalgebra. We adopt the procedure described in Sec. 3 of [11] to explore conjugacy classes.

*One-Dimensional Subalgebras*

From the description of the ideals, every one-dimensional subspace is of the form:

$$\langle X_{\alpha_1+\alpha_2}\rangle, \quad \langle X_{\alpha_2} + bX_{\alpha_1+\alpha_2}\rangle, \quad \text{and} \quad \langle X_{\alpha_1} + cX_{\alpha_2} + dX_{\alpha_1+\alpha_2}\rangle.$$

Since $\langle X_{\alpha_1+\alpha_2}\rangle$ is a one-dimensional ideal, it cannot be conjugate to any other subalgebra. Consider a one-dimensional subspace of the form $\langle X_{\alpha_1} + cX_{\alpha_2} + dX_{\alpha_1+\alpha_2}\rangle = \langle \overline{X_{\alpha_1}}\rangle$. The basis of $\mathfrak{g}'$ is $\overline{X_{\alpha_1}}, X_{\alpha_2}, X_{\alpha_1+\alpha_2}$. Since $\overline{X_{\alpha_1}}$ and $X_{\alpha_1} + aX_{\alpha_1+\alpha_2}$ centralize $\overline{X_{\alpha_1}}$, we consider the one-parameter subgroup generated by $X_{\alpha_2}$ to determine conjugates of $\langle \overline{X_{\alpha_1}}\rangle$. We have

$$ad(X_{\alpha_2})(\overline{X_{\alpha_1}}) = X_{\alpha_1+\alpha_2} \quad \text{and} \quad ad(X_{\alpha_2})(X_{\alpha_1+\alpha_2}) = 0.$$

Thus, $e^{\epsilon ad(X_{\alpha_2})}(\overline{X_{\alpha_1}}) = \overline{X_{\alpha_1}} + \epsilon X_{\alpha_1+\alpha_2}$. Setting $\epsilon = -d$ implies:

$$\langle X_{\alpha_1} + cX_{\alpha_2} + dX_{\alpha_1+\alpha_2}\rangle \sim \langle X_{\alpha_1} + cX_{\alpha_2}\rangle.$$

Similarly, $\langle X_{\alpha_2} + bX_{\alpha_1+\alpha_2}\rangle \sim \langle X_{\alpha_2}\rangle$.

Therefore, representatives of conjugacy classes of one-dimensional subalgebras are

$$\langle X_{\alpha2}\rangle, \quad \langle X_{\alpha1+\alpha2}\rangle, \quad \text{and} \quad \langle X_{\alpha1} + cX_{\alpha2}\rangle.$$

*Two-Dimensional Subalgebras*

A two-dimensional subalgebra is obtained by extending a one-dimensional subalgebra. For one-dimensional sbalgebra $\langle X_{\alpha1+\alpha2}\rangle$, the normalizer

$$N\langle X_{\alpha1+\alpha2}\rangle = \langle X_{\alpha1}, X_{\alpha2}, X_{\alpha1+\alpha2}\rangle.$$

Thus,

$$N\langle X_{\alpha1+\alpha2}\rangle/\langle X_{\alpha1+\alpha2}\rangle = \langle \overline{X_{\alpha1}}, \overline{X_{\alpha2}}\rangle,$$

which is abelian. Therefore, the one-dimensional subalgebras are $\langle \overline{X_{\alpha2}}\rangle$ and $\langle \overline{X_{\alpha1}} + \gamma\overline{X_{\alpha2}}\rangle$. The extensions of $\langle X_{\alpha1+\alpha2}\rangle$ are then $\langle X_{\alpha1+\alpha2}, X_{\alpha2}\rangle$ and $\langle X_{\alpha1+\alpha2}, X_{\alpha1} + \gamma X_{\alpha2}\rangle$.
For the one-dimensional subalgebra $\langle X_{\alpha2}\rangle$, we have

$$N\langle X_{\alpha2}\rangle/\langle X_{\alpha2}\rangle = \langle X_{\alpha1}, X_{\alpha1+\alpha2}\rangle/\langle X_{\alpha2}\rangle = \langle X_{\alpha1+\alpha2}\rangle.$$

Hence, $\langle X_{\alpha2}, X_{\alpha1+\alpha2}\rangle$ is the only extension of $\langle X_{\alpha2}\rangle$. Finally, for $\langle X_{\alpha_1} + cX_{\alpha_2}\rangle$, we find $N\langle X_{\alpha1} + cX_{\alpha2}\rangle = \langle X_{\alpha1} + cX_{\alpha2}, X_{\alpha1+\alpha2}\rangle$. Thus, the only extension is $\langle X_{\alpha1} + cX_{\alpha2}, X_{\alpha1+\alpha2}\rangle$.
Therefore, representatives of conjugacy classes of two-dimensional subalgebras of $A_2$ algebra are:

$$\langle X_{\alpha2}, X_{\alpha1+\alpha2}\rangle, \quad \text{and} \quad \langle X_{\alpha1+\alpha2}, X_{\alpha1} + \gamma X_{\alpha2}\rangle,$$

where $\gamma \in \mathbb{R}$.

### 4.2. *Borel subalgebras of type* $B_2$

For the Lie algebra $B_2$, corresponding to $so(5,\mathbb{C})$, the Maple code is executed with the root type `"B"` and $n$ set to $2$ as inputs. The code generates the Dynkin diagram for $B_2$, followed by the systematic construction of its positive roots, simple roots, and the Borel subalgebra.

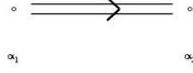

$$\text{Root\_Type} := B_2$$
$$\text{Positive\_Roots} := [\alpha 1, \alpha 2, \alpha 1 + \alpha 2, \alpha 1 + 2\alpha 2]$$
$$\text{Simple\_Roots} := [\alpha 1, \alpha 2]$$
$$\text{Derived\_Algebra} := [X_{\alpha 1}, X_{\alpha 2}, X_{\alpha 1 + \alpha 2}, X_{\alpha 1 + 2\alpha 2}]$$
$$\text{Borel\_Subalgebra} := [H_{\alpha 1}, H_{\alpha 2}, X_{\alpha 1}, X_{\alpha 2}, X_{\alpha 1 + \alpha 2}, X_{\alpha 1 + 2\alpha 2}]$$
$$\text{Ideals\_of\_DA\_of\_Borel\_Subalgebra} := \{[X_{\alpha 1 + 2\alpha 2}], [X_{\alpha 1 + \alpha 2}, X_{\alpha 1 + 2\alpha 2}], [X_{\alpha 1}, X_{\alpha 1 + \alpha 2}, X_{\alpha 1 + 2\alpha 2}],$$
$$[X_{\alpha 2}, X_{\alpha 1 + \alpha 2}, X_{\alpha 1 + 2\alpha 2}], [X_{\alpha 1}, X_{\alpha 2}, X_{\alpha 1 + \alpha 2}, X_{\alpha 1 + 2\alpha 2}]\}$$
$$\text{All\_Ideals\_of\_Borel\_Subalgebra} := \{[X_{\alpha 1 + 2\alpha 2}], [X_{\alpha 1 + \alpha 2}, X_{\alpha 1 + 2\alpha 2}], [X_{\alpha 1}, X_{\alpha 1 + \alpha 2}, X_{\alpha 1 + 2\alpha 2}],$$
$$[X_{\alpha 2}, X_{\alpha 1 + \alpha 2}, X_{\alpha 1 + 2\alpha 2}], [X_{\alpha 1}, X_{\alpha 2}, X_{\alpha 1 + \alpha 2}, X_{\alpha 1 + 2\alpha 2}], [H_{\alpha 2}, X_{\alpha 1}, X_{\alpha 2}, X_{\alpha 1 + \alpha 2}, X_{\alpha 1 + 2\alpha 2}],$$
$$[H_{\alpha 2} a_2 + H_{\alpha 1}, X_{\alpha 1}, X_{\alpha 2}, X_{\alpha 1 + \alpha 2}, X_{\alpha 1 + 2\alpha 2}], [H_{\alpha 1}, H_{\alpha 2}, X_{\alpha 1}, X_{\alpha 2}, X_{\alpha 1 + \alpha 2}, X_{\alpha 1 + 2\alpha 2}]\}$$

Figure 2: Ideals of a Borel subalgebras of type $B_2$

Moreover, the code identifies and lists all dimensional ideals of $B_2$ such as

$\langle X_{\alpha 1 + 2\alpha 2}\rangle$, $\langle X_{\alpha 1 + \alpha 2}, X_{\alpha 1 + 2\alpha 2}\rangle$, $\langle X_{\alpha 1}, X_{\alpha 1 + \alpha 2}, X_{\alpha 1 + 2\alpha 2}\rangle$, $\langle X_{\alpha 2}, X_{\alpha 1 + \alpha 2}, X_{\alpha 1 + 2\alpha 2}\rangle$, $\langle X_{\alpha 1}, X_{\alpha 2}, X_{\alpha 1 + \alpha 2}, X_{\alpha 1 + 2\alpha 2}\rangle$, $\langle H_{\alpha 2}, X_{\alpha 1}, X_{\alpha 2}, X_{\alpha 1 + \alpha 2}, X_{\alpha 1 + 2\alpha 2}\rangle$, $\langle H_{\alpha 1} + a_2 H_{\alpha 2}, X_{\alpha 1}, X_{\alpha 2}, X_{\alpha 1 + \alpha 2}, X_{\alpha 1 + 2\alpha 2}\rangle$, and $\langle H_{\alpha 1}, H_{\alpha 2}, X_{\alpha 1}, X_{\alpha 2}, X_{\alpha 1 + \alpha 2}, X_{\alpha 1 + 2\alpha 2}\rangle$.

The ideals of derived algebras are:

$$\langle X_{\alpha 1 + 2\alpha 2}\rangle, \quad \langle X_{\alpha 1 + \alpha 2}, X_{\alpha 1 + 2\alpha 2}\rangle, \quad \langle X_{\alpha 1}, X_{\alpha 1 + \alpha 2}, X_{\alpha 1 + 2\alpha 2}\rangle,$$
$$\langle X_{\alpha 2}, X_{\alpha 1 + \alpha 2}, X_{\alpha 1 + 2\alpha 2}\rangle, \quad \text{and} \quad \langle X_{\alpha 1}, X_{\alpha 2}, X_{\alpha 1 + \alpha 2}, X_{\alpha 1 + 2\alpha 2}\rangle.$$

The one-dimensional subspaces are given by:

$$\langle X_{\alpha 1 + 2\alpha 2}\rangle, \quad \langle X_{\alpha 1 + \alpha 2} + aX_{\alpha 1 + 2\alpha 2}\rangle, \quad \langle X_{\alpha 2} + dX_{\alpha 1 + \alpha 2} + eX_{\alpha 1 + 2\alpha 2}\rangle,$$
$$\text{and} \quad \langle X_{\alpha 1} + kX_{\alpha 2} + lX_{\alpha 1 + \alpha 2} + mX_{\alpha 1 + 2\alpha 2}\rangle.$$

*One-Dimensional Subalgebras*

To find one-dimensional subalgebras, we follow the same procedure as before. $\langle X_{\alpha 1 + 2\alpha 2}\rangle$ cannot jugate to any other subalgebras. In computing conjugates of $\langle X_{\alpha 1 + \alpha 2} + aX_{\alpha 1 + 2\alpha 2}\rangle$, we use the factorization of the adjoint group $e^{\langle X_{\alpha_2}\rangle}e^{\langle X_{\alpha_1+\alpha_2}+aX_{\alpha_1+2\alpha_2}\rangle}$. Now,

$$ad(X_{\alpha_2})(X_{\alpha_1+\alpha_2} + aX_{\alpha_1+2\alpha_2}) = X_{\alpha_1+2\alpha_2} \quad \text{and} \quad ad(X_{\alpha_2})(X_{\alpha_1+2\alpha_2}) = 0.$$

Thus,

$$e^{\epsilon ad(X_{\alpha_2})}(X_{\alpha_1+\alpha_2} + aX_{\alpha_1+2\alpha_2}) = X_{\alpha_1+\alpha_2} + aX_{\alpha_1+2\alpha_2} + \epsilon X_{\alpha_1+2\alpha_2}.$$

Let $\epsilon = -a$, which implies $\langle X_{\alpha_1+\alpha_2} + aX_{\alpha_1+2\alpha_2}\rangle \sim \langle X_{\alpha_1+\alpha_2}\rangle$.

By performing the same procedure, we obtain: $\langle X_{\alpha_2} + dX_{\alpha_1+\alpha_2} + eX_{\alpha_1+2\alpha_2}\rangle \sim \langle X_{\alpha_2}\rangle$, and $\langle X_{\alpha_1} + kX_{\alpha_2} + lX_{\alpha_1+\alpha_2} + mX_{\alpha_1+2\alpha_2}\rangle \sim \langle X_{\alpha_1} + kX_{\alpha_2}\rangle$.

Therefore, representatives of conjugacy classes of one-dimensional subalgebras are

$$\langle X_{\alpha_1+2\alpha_2}\rangle, \quad \langle X_{\alpha_1+\alpha_2}\rangle, \quad \langle X_{\alpha_2}\rangle, \quad \text{and} \quad \langle X_{\alpha_1} + kX_{\alpha_2}\rangle.$$

*Two-Dimensional Subalgebras:*

For the one-dimensional subalgebra $\langle X_{\alpha_2}\rangle$,

$$N\langle X_{\alpha_2}\rangle = \langle X_{\alpha_2}, X_{\alpha_1+2\alpha_2}\rangle.$$

The quotient

$$N\langle X_{\alpha_2}\rangle/\langle X_{\alpha_2}\rangle = \langle X_{\alpha_1+2\alpha_2}\rangle.$$

Thus, the only extension is $\langle X_{\alpha_1+2\alpha_2}, X_{\alpha_2}\rangle$.

For the one-dimensional subalgebra $\langle X_{\alpha_1+\alpha_2}\rangle$.

$$N\langle X_{\alpha_1+\alpha_2}\rangle = \langle X_{\alpha_1}, X_{\alpha_1+\alpha_2}, X_{\alpha_1+2\alpha_2}\rangle$$

The quotient

$$N\langle X_{\alpha_1+\alpha_2}\rangle/\langle X_{\alpha_1+\alpha_2}\rangle = \langle X_{\alpha_1}, X_{\alpha_1+2\alpha_2}\rangle,$$

is abelian. Thus, the extensions of $\langle X_{\alpha_1+\alpha_2}\rangle$ are

$$\langle X_{\alpha_1+\alpha_2}, X_{\alpha_1+2\alpha_2}\rangle \quad \text{and} \quad \langle X_{\alpha_1+\alpha_2}, X_{\alpha_1} + \lambda X_{\alpha_1+2\alpha_2}\rangle.$$

For $N\langle X_{\alpha_1} + kX_{\alpha_2}\rangle/\langle X_{\alpha_1} + kX_{\alpha_2}\rangle = \langle X_{\alpha_1+2\alpha_2}\rangle$, the only extension is $\langle X_{\alpha_1+2\alpha_2}, X_{\alpha_1} + kX_{\alpha_2}\rangle$.

The quotient $N\langle X_{\alpha_1+2\alpha_2}\rangle/\langle X_{\alpha_1+2\alpha_2}\rangle = \langle X_{\alpha_1}, X_{\alpha_2}, X_{\alpha_1+\alpha_2}\rangle$ is a three-dimensional non abelian al-gebra $[X_{\alpha1}, X_{\alpha2}] = X_{\alpha1+\alpha2}$, implies the one-dimensional representations are $\langle X_{\alpha1+\alpha2}\rangle$, $\langle X_{\alpha2}\rangle$ and $\langle X_{\alpha1} + \mu X_{\alpha2}\rangle$. Thus the extensions are $\langle X_{\alpha1+\alpha2}, X_{\alpha1+2\alpha2}\rangle$, $\langle X_{\alpha1+2\alpha2}, X_{\alpha2}\rangle$ and $\langle X_{\alpha1+2\alpha2}, X_{\alpha1} + \mu X_{\alpha2}\rangle$.

Therefore the two dimensional subalgebras are

$$\langle X_{\alpha1+2\alpha2}, X_{\alpha2}\rangle, \quad \langle X_{\alpha1+\alpha2}, X_{\alpha1+2\alpha2}\rangle, \quad \langle X_{\alpha1+\alpha2}, X_{\alpha1} + \lambda X_{\alpha1+2\alpha2}\rangle, \quad \text{and} \quad \langle X_{\alpha1+2\alpha2}, X_{\alpha1} + \mu X_{\alpha2}\rangle.$$

*Three-Dimensional Subalgebras*

To obtain three-dimensional subalgebras, we need to compute the one-dimensional extensions of representatives from the conjugacy classes of two-dimensional subalgebras. By performing this for all classes of two-dimensional algebras, we obtain the following extensions:

For two-dimensional subalgebra $\langle X_{\alpha1+2\alpha2}, X_{\alpha2}\rangle$,

$$N\langle X_{\alpha1+2\alpha2}, X_{\alpha2}\rangle = \langle X_{\alpha1+2\alpha2}, X_{\alpha1+\alpha2}, X_{\alpha2}\rangle.$$

The quotient

$$N\langle X_{\alpha2}, X_{\alpha1+2\alpha2}\rangle / \langle X_{\alpha2}, X_{\alpha1+2\alpha2}\rangle = \langle X_{\alpha1+\alpha2}\rangle,$$

implying the only extension is

$$\langle X_{\alpha2}, X_{\alpha1+\alpha2}, X_{\alpha1+2\alpha2}\rangle.$$

Similarly, applying the same procedure to other classes of two-dimensional algebras, we obtain the following three-dimensional subalgebra:

$$\langle X_{\alpha1} + \mu X_{\alpha2}, X_{\alpha1+\alpha2}, X_{\alpha1+2\alpha2}\rangle.$$

Therefore, the representatives of conjugacy classes of three-dimensional subalgebras of derived algebra of $B_2$ algebra are

$$\langle X_{\alpha2}, X_{\alpha1+\alpha2}, X_{\alpha1+2\alpha2}\rangle, \quad \text{and} \quad \langle X_{\alpha1} + \mu X_{\alpha2}, X_{\alpha1+\alpha2}, X_{\alpha1+2\alpha2}\rangle.$$

*4.3. Borel subalgebras of type $G_2$*

For the exceptional Lie algebra $G_2$, the code is executed with root type "G" and $n$ is 2. The Dynkin diagram, representing the structure of $G_2$, is plotted and the positive roots, simple roots, and Borel subalgebra are constructed systematically. $G_2$ is irreducible and has a unique one-dimensional ideal in its Borel subalgebra. The code then determines all the higher-dimensional ideals, providing a complete classification of the ideals of $G_2$. The results are displayed in Figure 3.

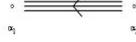

$$Root\_Type := G_2$$
$$Positive\_Roots := [\alpha 1, \alpha 2, \alpha 1 + \alpha 2, 2\alpha 1 + \alpha 2, 3\alpha 1 + \alpha 2, 3\alpha 1 + 2\alpha 2]$$
$$Simple\_Roots := [\alpha 1, \alpha 2]$$
$$Derived\_Algebra := [X_{\alpha 1}, X_{\alpha 2}, X_{\alpha 1+\alpha 2}, X_{2\alpha 1+\alpha 2}, X_{3\alpha 1+\alpha 2}, X_{3\alpha 1+2\alpha 2}]$$
$$Borel\_Subalgebra := [H_{\alpha 1}, H_{\alpha 2}, X_{\alpha 1}, X_{\alpha 2}, X_{\alpha 1+\alpha 2}, X_{2\alpha 1+\alpha 2}, X_{3\alpha 1+\alpha 2}, X_{3\alpha 1+2\alpha 2}]$$
$$Ideals\_of\_DA\_of\_Borel\_Subalgebra := \{[X_{3\alpha 1+2\alpha 2}], [X_{3\alpha 1+\alpha 2}, X_{3\alpha 1+2\alpha 2}], [X_{2\alpha 1+\alpha 2}, X_{3\alpha 1+\alpha 2}, X_{3\alpha 1+2\alpha 2}], [X_{\alpha 1+\alpha 2}, X_{2\alpha 1+\alpha 2}, X_{3\alpha 1+\alpha 2},$$
$$X_{3\alpha 1+2\alpha 2}], [X_{\alpha 1}, X_{\alpha 1+\alpha 2}, X_{2\alpha 1+\alpha 2}, X_{3\alpha 1+\alpha 2}, X_{3\alpha 1+2\alpha 2}], [X_{\alpha 2}, X_{\alpha 1+\alpha 2}, X_{2\alpha 1+\alpha 2}, X_{3\alpha 1+\alpha 2}, X_{3\alpha 1+2\alpha 2}], [X_{\alpha 1}, X_{\alpha 2}, X_{\alpha 1+\alpha 2}, X_{2\alpha 1+\alpha 2},$$
$$X_{3\alpha 1+\alpha 2}, X_{3\alpha 1+2\alpha 2}]\}$$
$$All\_Ideals\_of\_Borel\_Subalgebra := \{[X_{3\alpha 1+2\alpha 2}], [X_{3\alpha 1+\alpha 2}, X_{3\alpha 1+2\alpha 2}], [X_{2\alpha 1+\alpha 2}, X_{3\alpha 1+\alpha 2}, X_{3\alpha 1+2\alpha 2}], [X_{\alpha 1+\alpha 2}, X_{2\alpha 1+\alpha 2}, X_{3\alpha 1+\alpha 2},$$
$$X_{3\alpha 1+2\alpha 2}], [X_{\alpha 1}, X_{\alpha 1+\alpha 2}, X_{2\alpha 1+\alpha 2}, X_{3\alpha 1+\alpha 2}, X_{3\alpha 1+2\alpha 2}], [X_{\alpha 2}, X_{\alpha 1+\alpha 2}, X_{2\alpha 1+\alpha 2}, X_{3\alpha 1+\alpha 2}, X_{3\alpha 1+2\alpha 2}], [X_{\alpha 1}, X_{\alpha 2}, X_{\alpha 1+\alpha 2}, X_{2\alpha 1+\alpha 2},$$
$$X_{3\alpha 1+\alpha 2}, X_{3\alpha 1+2\alpha 2}], [H_{\alpha 2}, X_{\alpha 1}, X_{\alpha 2}, X_{\alpha 1+\alpha 2}, X_{2\alpha 1+\alpha 2}, X_{3\alpha 1+\alpha 2}, X_{3\alpha 1+2\alpha 2}], [H_{\alpha 2}a_2 + H_{\alpha 1}, X_{\alpha 1}, X_{\alpha 2}, X_{\alpha 1+\alpha 2}, X_{2\alpha 1+\alpha 2}, X_{3\alpha 1+\alpha 2},$$
$$X_{3\alpha 1+2\alpha 2}], [H_{\alpha 1}, H_{\alpha 2}, X_{\alpha 1}, X_{\alpha 2}, X_{\alpha 1+\alpha 2}, X_{2\alpha 1+\alpha 2}, X_{3\alpha 1+\alpha 2}, X_{3\alpha 1+2\alpha 2}]\}$$

Figure 3: Ideals of a Borel subalgebra of type $G_2$

Using the ideals of derived algebras of the Borel subalgebra of the $G_2$ algebra. Representatives of conjugacy classes of one-dimensional subalgebras are

$$\langle X_{3\alpha 1+2\alpha 2}\rangle, \quad \langle X_{3\alpha 1+\alpha 2}\rangle, \quad \langle X_{2\alpha 1+\alpha 2}\rangle, \quad \langle X_{\alpha 1+\alpha 2}\rangle, \quad \langle X_{\alpha 2}\rangle, \quad \text{and} \quad \langle X_{\alpha 1}+\beta X_{\alpha 2}\rangle.$$

Representatives of conjugacy classes of two-dimensional subalgebras are

$$\langle X_{\alpha 1}+AX_{\alpha 2}, X_{3\alpha 1+2\alpha 2}\rangle, \quad \langle X_{\alpha 2}, X_{3\alpha 1+2\alpha 2}\rangle, \quad \langle X_{\alpha 2}, X_{2\alpha 1+\alpha 2}\rangle, \quad \langle X_{\alpha 2}, X_{\alpha 1+\alpha 2}+BX_{2\alpha 1+\alpha 2}\rangle,$$
$$\langle X_{\alpha 1+\alpha 2}, X_{3\alpha 1+2\alpha 2}\rangle, \quad \langle X_{\alpha 1+\alpha 2}, X_{3\alpha 1+\alpha 2}\rangle, \quad \langle X_{\alpha 1+\alpha 2}, X_{\alpha 2}+DX_{3\alpha 1+\alpha 2}\rangle, \quad \langle X_{2\alpha 1+\alpha 2}, X_{3\alpha 1+2\alpha 2}\rangle,$$
$$\langle X_{2\alpha 1+\alpha 2}, X_{3\alpha 1+\alpha 2}\rangle, \quad \langle X_{2\alpha 1+\alpha 2}, X_{\alpha 2}+EX_{3\alpha 1+\alpha 2}\rangle, \quad \langle X_{3\alpha 1+\alpha 2}, X_{3\alpha 1+2\alpha 2}\rangle, \quad \text{and} \quad \langle X_{\alpha 1}+FX_{\alpha 1+\alpha 2}, X_{3\alpha 1+\alpha 2}\rangle.$$

Representatives of conjugacy classes of three-dimensional subalgebras are $\langle X_{\alpha 1} + b_1 X_{\alpha 2}, X_{3\alpha 1+\alpha 2}, X_{3\alpha 1+2\alpha 2} \rangle$, $\langle X_{\alpha 2}, X_{3\alpha 1+\alpha 2}, X_{3\alpha 1+2\alpha 2} \rangle$, $\langle X_{2\alpha 1+\alpha 2} + b_2 X_{3\alpha 1+\alpha 2}, X_{\alpha 2}, X_{3\alpha 1+2\alpha 2} \rangle$, $\langle X_{\alpha 1+\alpha 2} + b_3 X_{2\alpha 1+\alpha 2} + b_4 X_{3\alpha 1+\alpha 2}, X_{\alpha 2}, X_{3\alpha 1+2\alpha 2} \rangle$, $\langle X_{\alpha 2}, X_{2\alpha 1+\alpha 2}, X_{3\alpha 1+2\alpha 2} \rangle$, $\langle X_{\alpha 1+\alpha 2} + b_5 X_{2\alpha 1+\alpha 2}, X_{\alpha 2}, X_{3\alpha 1+2\alpha 2} \rangle$, $\langle X_{\alpha 1+\alpha 2}, X_{3\alpha 1+\alpha 2}, X_{3\alpha 1+2\alpha 2} \rangle$, $\langle X_{2\alpha 1+\alpha 2} + b_6 X_{3\alpha 1+\alpha 2}, X_{\alpha 1+\alpha 2}, X_{3\alpha 1+2\alpha 2} \rangle$, $\langle X_{\alpha 2} + b_7 X_{2\alpha 1+\alpha 2} + b_8 X_{3\alpha 1+\alpha 2}, X_{\alpha 1+\alpha 2}, X_{3\alpha 1+2\alpha 2} \rangle$, $\langle X_{\alpha 2} + b_9 X_{3\alpha 1+\alpha 2}, X_{\alpha 1+\alpha 2}, X_{3\alpha 1+2\alpha 2} \rangle$, $\langle X_{2\alpha 1+\alpha 2}, X_{3\alpha 1+\alpha 2}, X_{3\alpha 1+2\alpha 2} \rangle$, $\langle X_{\alpha 1+\alpha 2} + b_{10} X_{3\alpha 1+\alpha 2}, X_{2\alpha 1+\alpha 2}, X_{3\alpha 1+2\alpha 2} \rangle$, $\langle X_{\alpha 2} + b_{11} X_{\alpha 1+\alpha 2} + b_{12} X_{3\alpha 1+\alpha 2}, X_{2\alpha 1+\alpha 2}, X_{3\alpha 1+2\alpha 2} \rangle$, $\langle X_{\alpha 1} + b_{13} X_{3\alpha 1+2\alpha 2}, X_{2\alpha 1+\alpha 2}, X_{3\alpha 1+\alpha 2} \rangle$, $\langle X_{\alpha 2} + b_{14} X_{3\alpha 1+\alpha 2}, X_{2\alpha 1+\alpha 2}, X_{3\alpha 1+2\alpha 2} \rangle$, and $\langle X_{\alpha 1} + b_{15} X_{\alpha 1+\alpha 2}, X_{2\alpha 1+\alpha 2}, X_{3\alpha 1+2\alpha 2} \rangle$.

## 5. Abelian Ideals of type $F_4$ Algebra

Using the code described in Sec. 3, we determined that the Borel subalgebra of type $F_4$ contains a total of 110 ideals. The code executed in approximately 2 seconds for this computation. Among these, 16 are abelian ideals. Below, we present the complete list of all abelian ideals in the Borel subalgebra of type $F_4$.

1. $0$
2. $\langle X_{2\alpha_1+3\alpha_2+4\alpha_3+2\alpha_4} \rangle$
3. $\langle X_{\alpha_1+3\alpha_2+4\alpha_3+2\alpha_4}, X_{2\alpha_1+3\alpha_2+4\alpha_3+2\alpha_4} \rangle$
4. $\langle X_{\alpha_1+2\alpha_2+4\alpha_3+2\alpha_4}, X_{\alpha_1+3\alpha_2+4\alpha_3+2\alpha_4}, X_{2\alpha_1+3\alpha_2+4\alpha_3+2\alpha_4} \rangle$
5. $\langle X_{\alpha_1+2\alpha_2+3\alpha_3+2\alpha_4}, X_{\alpha_1+2\alpha_2+4\alpha_3+2\alpha_4}, X_{\alpha_1+3\alpha_2+4\alpha_3+2\alpha_4}, X_{2\alpha_1+3\alpha_2+4\alpha_3+2\alpha_4} \rangle$
6. $\langle X_{\alpha_1+2\alpha_2+2\alpha_3+2\alpha_4}, X_{\alpha_1+2\alpha_2+3\alpha_3+2\alpha_4}, X_{\alpha_1+2\alpha_2+4\alpha_3+2\alpha_4}, X_{\alpha_1+3\alpha_2+4\alpha_3+2\alpha_4},$
   $X_{2\alpha_1+3\alpha_2+4\alpha_3+2\alpha_4} \rangle$
7. $\langle X_{\alpha_1+2\alpha_2+3\alpha_3+\alpha_4}, X_{\alpha_1+2\alpha_2+3\alpha_3+2\alpha_4}, X_{\alpha_1+2\alpha_2+4\alpha_3+2\alpha_4}, X_{\alpha_1+3\alpha_2+4\alpha_3+2\alpha_4},$
   $X_{2\alpha_1+3\alpha_2+4\alpha_3+2\alpha_4} \rangle$
8. $\langle X_{\alpha_1+\alpha_2+2\alpha_3+2\alpha_4}, X_{\alpha_1+2\alpha_2+2\alpha_3+2\alpha_4}, X_{\alpha_1+2\alpha_2+3\alpha_3+2\alpha_4}, X_{\alpha_1+2\alpha_2+4\alpha_3+2\alpha_4},$
   $X_{\alpha_1+3\alpha_2+4\alpha_3+2\alpha_4}, X_{2\alpha_1+3\alpha_2+4\alpha_3+2\alpha_4} \rangle$
9. $\langle X_{\alpha_1+2\alpha_2+2\alpha_3+2\alpha_4}, X_{\alpha_1+2\alpha_2+3\alpha_3+\alpha_4}, X_{\alpha_1+2\alpha_2+3\alpha_3+2\alpha_4}, X_{\alpha_1+2\alpha_2+4\alpha_3+2\alpha_4},$
   $X_{\alpha_1+3\alpha_2+4\alpha_3+2\alpha_4}, X_{2\alpha_1+3\alpha_2+4\alpha_3+2\alpha_4} \rangle$

10. $\langle X_{\alpha_2+2\alpha_3+2\alpha_4}, X_{\alpha_1+\alpha_2+2\alpha_3+2\alpha_4}, X_{\alpha_1+2\alpha_2+2\alpha_3+2\alpha_4}, X_{\alpha_1+2\alpha_2+3\alpha_3+2\alpha_4},$
$X_{\alpha_1+2\alpha_2+4\alpha_3+2\alpha_4}, X_{\alpha_1+3\alpha_2+4\alpha_3+2\alpha_4}, X_{2\alpha_1+3\alpha_2+4\alpha_3+2\alpha_4}\rangle$

11. $\langle X_{\alpha1+\alpha2+2\alpha3+2\alpha4}, X_{\alpha1+2\alpha2+2\alpha3+2\alpha4}, X_{\alpha1+2\alpha2+3\alpha3+\alpha4}, X_{\alpha1+2\alpha2+3\alpha3+2\alpha4},$
$X_{\alpha1+2\alpha2+4\alpha3+2\alpha4}, X_{\alpha1+3\alpha2+4\alpha3+2\alpha4}, X_{2\alpha1+3\alpha2+4\alpha3+2\alpha4}\rangle$

12. $\langle X_{\alpha1+2\alpha2+2\alpha3+\alpha4}, X_{\alpha1+2\alpha2+2\alpha3+2\alpha4}, X_{\alpha1+2\alpha2+3\alpha3+\alpha4}, X_{\alpha1+2\alpha2+3\alpha3+2\alpha4},$
$X_{\alpha1+2\alpha2+4\alpha3+2\alpha4}, X_{\alpha1+3\alpha2+4\alpha3+2\alpha4}, X_{2\alpha1+3\alpha2+4\alpha3+2\alpha4}\rangle$

13. $\langle X_{\alpha1+2\alpha2+2\alpha3}, X_{\alpha1+2\alpha2+2\alpha3+\alpha4}, X_{\alpha1+2\alpha2+2\alpha3+2\alpha4}, X_{\alpha1+2\alpha2+3\alpha3+\alpha4},$
$X_{\alpha1+2\alpha2+3\alpha3+2\alpha4}, X_{\alpha1+2\alpha2+4\alpha3+2\alpha4}, X_{\alpha1+3\alpha2+4\alpha3+2\alpha4}, X_{2\alpha1+3\alpha2+4\alpha3+2\alpha4}\rangle$

14. $\langle X_{\alpha2+2\alpha3+2\alpha4}, X_{\alpha1+\alpha2+2\alpha3+2\alpha4}, X_{\alpha1+2\alpha2+2\alpha3+2\alpha4}, X_{\alpha1+2\alpha2+3\alpha3+\alpha4},$
$X_{\alpha1+2\alpha2+3\alpha3+2\alpha4}, X_{\alpha1+2\alpha2+4\alpha3+2\alpha4}, X_{\alpha1+3\alpha2+4\alpha3+2\alpha4}, X_{2\alpha1+3\alpha2+4\alpha3+2\alpha4}\rangle$

15. $\langle X_{\alpha1+\alpha2+2\alpha3+2\alpha4}, X_{\alpha1+2\alpha2+2\alpha3+\alpha4}, X_{\alpha1+2\alpha2+2\alpha3+2\alpha4}, X_{\alpha1+2\alpha2+3\alpha3+\alpha4},$
$X_{\alpha1+2\alpha2+3\alpha3+2\alpha4}, X_{\alpha1+2\alpha2+4\alpha3+2\alpha4}, X_{\alpha1+3\alpha2+4\alpha3+2\alpha4}, X_{2\alpha1+3\alpha2+4\alpha3+2\alpha4}\rangle$

16. $\langle X_{\alpha2+2\alpha3+2\alpha4}, X_{\alpha1+\alpha2+2\alpha3+2\alpha4}, X_{\alpha1+2\alpha2+2\alpha3+\alpha4}, X_{\alpha1+2\alpha2+2\alpha3+2\alpha4},$
$X_{\alpha1+2\alpha2+3\alpha3+\alpha4}, X_{\alpha1+2\alpha2+3\alpha3+2\alpha4}, X_{\alpha1+2\alpha2+4\alpha3+2\alpha4}, X_{\alpha1+3\alpha2+4\alpha3+2\alpha4},$
$X_{2\alpha1+3\alpha2+4\alpha3+2\alpha4}\rangle$

## Appendix

Below are the functions used in the above-mentioned code:

```
ConvertToList := proc (v)
    if type(v, 'Vector') then
        return [seq(v[i], i = 1 .. LinearAlgebra[Dimension](v))]
    else
        return v
    end if;
end proc:
```

```
AlphaForge := proc (proot, sroot)
    local dim, alphanames, sroottoalpha, Simple_root, Positive_root, i, j,
    coefficients, linear_combination;
    dim := LinearAlgebra[Dimension](sroot[1]);
```

```
4       alphanames := [seq(cat(alpha, i), i = 1 .. nops(sroot))];
5       sroottoalpha := table();
6       for i to nops(sroot) do
7           sroottoalpha[ConvertToList(sroot[i])] := parse(alphanames[i]) end do;
8           Simple_root := [seq(sroottoalpha[ConvertToList(sroot[i])], i = 1 .. nops(
        sroot))];
9           Positive_root := [];
10          for i to nops(proot) do
11              coefficients := ConvertToList(proot[i]);
12              linear_combination := 0;
13              for j to nops(sroot) do
14                  if coefficients[j] <> 0 then
15                      linear_combination := linear_combination+coefficients[j]*
        sroottoalpha[ConvertToList(sroot[j])]
16                  end if
17              end do;
18              Positive_root := [op(Positive_root), simplify(linear_combination)];
19          end do;
20          return Positive_root, Simple_root;
21      end proc:
```

```
1   ConstructCartanSubalgebra := proc (Positive_Roots, Simple_Roots)
2       local Cartan_Subalgebra, i;
3           Cartan_Subalgebra := [seq(cat(H[Simple_Roots[i]]), i = 1 .. nops(
        Simple_Roots))];
4       return Cartan_Subalgebra;
5   end proc:
```

```
1   ConstructDerivedAlgebra := proc (Positive_Roots, Simple_Roots)
2       local DerivedAlgebra, i;
3           DerivedAlgebra := [seq(cat(X[Positive_Roots[i]]), i = 1 .. nops(
        Positive_Roots))];
4       return DerivedAlgebra;
5   end proc:
```

```
1   ConstructBorelSubalgebra := proc (CS, DA)
2       local BorelSubalgebra;
3           BorelSubalgebra := [op(CS), op(DA)];
4           return BorelSubalgebra;
5   end proc:
```

```
1   lieB := proc (Xalpha, Xbeta, proot)
2       local alpha, beta, sum_of_roots;
```

```
3            alpha := op(1, Xalpha);
4            beta := op(1, Xbeta);
5            sum_of_roots := alpha+beta;
6            if `in`(sum_of_roots, proot) then
7                return X[sum_of_roots]
8            elif sum_of_roots = 0 or not `in`(sum_of_roots, proot) then
9                return 0
10           else error "Invalid input or root configuration."
11           end if;
12 end proc:
```

```
1 sortList := proc (sublist)
2     sort(sublist);
3 end proc:
```

```
1 CreateBasisList := proc(inputList::list)
2   local n, result;
3     n := nops(inputList);
4     result := [seq(inputList[i .. n], i = n .. 1, -1)];
5   return result
6 end proc:
```

```
1 BasisListCartan := proc (inputList::list)
2     local n, result; n := nops(inputList);
3     result := [seq(inputList[i .. n], i = n .. 1, -1)];
4     return result;
5 end proc:
```

```
1 CoefficientListCartan := proc (originalList)
2     local resultList, i;
3     resultList := [[1]];
4     for i to nops(originalList)-1 do
5         resultList := [op(resultList), [1, op(originalList[-i .. -1])]]
6     end do;
7     return resultList;
8 end proc:
```

```
1 DotProductList := proc (u::list, v::list)
2     local i; return [add(u[i]*v[i], i = 1 .. min(nops(u), nops(v)))];
3 end proc:
```

```
1 CalculateDotProducts := proc (cefflist::(list(list)), basislist::(list(list)))
```

```
2        return zip(proc (a, b) options operator, arrow; DotProductList(a, b) end proc,
    cefflist, basislist)
3 end proc:
```

```
1 DetermineOneDimensionalIdeal := proc (proot, sroot)
2     local i, j, IsZeroForAll, result, OneDimensionalIdeals;
3     OneDimensionalIdeals := [];
4     for i in proot do
5         IsZeroForAll := true;
6         for j in sroot do
7             result := lieB(X[j], X[i], proot);
8             if result <> 0 then
9                 IsZeroForAll := false;
10                break;
11            end if;
12        end do;
13        if IsZeroForAll then
14            OneDimensionalIdeals := {op(OneDimensionalIdeals), [X[i]]};
15        end if;
16    end do;
17    return OneDimensionalIdeals;
18 end proc:
```

```
1 RootsExtractor := proc (ideal)
2     local roots, m;
3     roots := [];
4     for m in ideal do
5         roots := [op(roots), op(1, m)];
6     end do;
7     return roots;
8 end proc:
```

```
1 DetermineOtherIdeals := proc (I_deal, proot, sroot)
2    local k, l, IsZeroForAll, result, nextIdeal;
3    nextIdeal := [];
4    for k in proot do
5        IsZeroForAll := true;
6        for l in sroot do
7            result := lieB(X[l], X[k], proot);
8            if result <> 0 then
9                IsZeroForAll := false;
10               break;
```

```
11            end if;
12        end do;
13        if IsZeroForAll then
14            nextIdeal := {op(nextIdeal), [X[k], op(I_deal)]};
15        end if;
16    end do;
17 return nextIdeal;
18 end proc:
```